\documentclass[12pt]{article}
\usepackage{amsmath,amssymb,amsbsy,amsfonts,
latexsym,amsopn,amstext,amsxtra}

\begin{document}

\newtheorem{theorem}{Theorem}
\newtheorem{lemma}[theorem]{Lemma}
\newtheorem{claim}[theorem]{Claim}
\newtheorem{cor}[theorem]{Corollary}
\newtheorem{prop}[theorem]{Proposition}
\newtheorem{definition}{Definition}
\newtheorem{question}[theorem]{Question}

\newenvironment{proof}{\noindent{\it Proof. }}{{\qed}}

\def\squareforqed{\hbox{\rlap{$\sqcap$}$\sqcup$}}
\def\qed{\ifmmode\squareforqed\else{\unskip\nobreak\hfil
\penalty50\hskip1em\null\nobreak\hfil\squareforqed
\parfillskip=0pt\finalhyphendemerits=0\endgraf}\fi}

\newfont{\teneufm}{eufm10}
\newfont{\seveneufm}{eufm7}
\newfont{\fiveeufm}{eufm5}
%
%
\newfam\eufmfam
     \textfont\eufmfam=\teneufm
\scriptfont\eufmfam=\seveneufm
     \scriptscriptfont\eufmfam=\fiveeufm
%
%
\def\frak#1{{\fam\eufmfam\relax#1}}

\def\eqref#1{(\ref{#1})}


\def\bbbr{{\rm I\!R}} 
\def\bbbm{{\rm I\!M}}
\def\bbbn{{\rm I\!N}} 
\def\bbbf{{\rm I\!F}}
\def\bbbh{{\rm I\!H}}
\def\bbbk{{\rm I\!K}}
\def\bbbp{{\rm I\!P}}

\def\vec#1{\mathbf{#1}}

\def\bbbone{{\mathchoice {\rm 1\mskip-4mu l} {\rm
1\mskip-4mu l}
{\rm 1\mskip-4.5mu l} {\rm 1\mskip-5mu l}}}
\def\bbbc{{\mathchoice
{\setbox0=\hbox{$\displaystyle\rm C$}\hbox{\hbox
to0pt{\kern0.4\wd0\vrule height0.9\ht0\hss}\box0}}
{\setbox0=\hbox{$\textstyle\rm C$}\hbox{\hbox
to0pt{\kern0.4\wd0\vrule height0.9\ht0\hss}\box0}}
{\setbox0=\hbox{$\scriptstyle\rm C$}\hbox{\hbox
to0pt{\kern0.4\wd0\vrule height0.9\ht0\hss}\box0}}
{\setbox0=\hbox{$\scriptscriptstyle\rm C$}\hbox{\hbox
to0pt{\kern0.4\wd0\vrule height0.9\ht0\hss}\box0}}}}
\def\bbbq{{\mathchoice
{\setbox0=\hbox{$\displaystyle\rm
Q$}\hbox{\raise
0.15\ht0\hbox to0pt{\kern0.4\wd0\vrule
height0.8\ht0\hss}\box0}}
{\setbox0=\hbox{$\textstyle\rm Q$}\hbox{\raise
0.15\ht0\hbox to0pt{\kern0.4\wd0\vrule
height0.8\ht0\hss}\box0}}
{\setbox0=\hbox{$\scriptstyle\rm Q$}\hbox{\raise
0.15\ht0\hbox to0pt{\kern0.4\wd0\vrule
height0.7\ht0\hss}\box0}}
{\setbox0=\hbox{$\scriptscriptstyle\rm Q$}\hbox{\raise
0.15\ht0\hbox to0pt{\kern0.4\wd0\vrule
height0.7\ht0\hss}\box0}}}}
\def\bbbt{{\mathchoice
{\setbox0=\hbox{$\displaystyle\rm
T$}\hbox{\hbox to0pt{\kern0.3\wd0\vrule
height0.9\ht0\hss}\box0}}
{\setbox0=\hbox{$\textstyle\rm T$}\hbox{\hbox
to0pt{\kern0.3\wd0\vrule height0.9\ht0\hss}\box0}}
{\setbox0=\hbox{$\scriptstyle\rm T$}\hbox{\hbox
to0pt{\kern0.3\wd0\vrule height0.9\ht0\hss}\box0}}
{\setbox0=\hbox{$\scriptscriptstyle\rm T$}\hbox{\hbox
to0pt{\kern0.3\wd0\vrule height0.9\ht0\hss}\box0}}}}
\def\bbbs{{\mathchoice
{\setbox0=\hbox{$\displaystyle     \rm
S$}\hbox{\raise0.5\ht0\hbox
to0pt{\kern0.35\wd0\vrule height0.45\ht0\hss}\hbox
to0pt{\kern0.55\wd0\vrule height0.5\ht0\hss}\box0}}
{\setbox0=\hbox{$\textstyle        \rm
S$}\hbox{\raise0.5\ht0\hbox
to0pt{\kern0.35\wd0\vrule height0.45\ht0\hss}\hbox
to0pt{\kern0.55\wd0\vrule height0.5\ht0\hss}\box0}}
{\setbox0=\hbox{$\scriptstyle      \rm
S$}\hbox{\raise0.5\ht0\hbox
to0pt{\kern0.35\wd0\vrule
height0.45\ht0\hss}\raise0.05\ht0\hbox
to0pt{\kern0.5\wd0\vrule height0.45\ht0\hss}\box0}}
{\setbox0=\hbox{$\scriptscriptstyle\rm
S$}\hbox{\raise0.5\ht0\hbox
to0pt{\kern0.4\wd0\vrule
height0.45\ht0\hss}\raise0.05\ht0\hbox
to0pt{\kern0.55\wd0\vrule height0.45\ht0\hss}\box0}}}}
\def\bbbz{{\mathchoice {\hbox{$\sf\textstyle
Z\kern-0.4em Z$}}
{\hbox{$\sf\textstyle Z\kern-0.4em Z$}}
{\hbox{$\sf\scriptstyle Z\kern-0.3em Z$}}
{\hbox{$\sf\scriptscriptstyle Z\kern-0.2em Z$}}}}
\def\ts{\thinspace}

\def\squareforqed{\hbox{\rlap{$\sqcap$}$\sqcup$}}
\def\qed{\ifmmode\squareforqed\else{\unskip\nobreak\hfil
\penalty50\hskip1em\null\nobreak\hfil\squareforqed
\parfillskip=0pt\finalhyphendemerits=0\endgraf}\fi}

\def\cA{{\mathcal A}}
\def\cB{{\mathcal B}}
\def\cC{{\mathcal C}}
\def\cD{{\mathcal D}}
\def\cE{{\mathcal E}}
\def\cF{{\mathcal F}}
\def\cG{{\mathcal G}}
\def\cH{{\mathcal H}}
\def\cI{{\mathcal I}}
\def\cJ{{\mathcal J}}
\def\cK{{\mathcal K}}
\def\cL{{\mathcal L}}
\def\cM{{\mathcal M}}
\def\cN{{\mathcal N}}
\def\cO{{\mathcal O}}
\def\cP{{\mathcal P}}
\def\cQ{{\mathcal Q}}
\def\cR{{\mathcal R}}
\def\cS{{\mathcal S}}
\def\cT{{\mathcal T}}
\def\cU{{\mathcal U}}
\def\cV{{\mathcal V}}
\def\cW{{\mathcal W}}
\def\cX{{\mathcal X}}
\def\cY{{\mathcal Y}}
\def\cZ{{\mathcal Z}}
\newcommand{\rmod}[1]{\: \mbox{mod} \: #1}

\def\vr{\mathbf r}

\def\e{{\mathbf{\,e}}}
\def\ep{{\mathbf{\,e}}_p}
\def\em{{\mathbf{\,e}}_m}

\def\Tr{{\mathrm{Tr}}}
\def\ind{{\mathrm{ind}\,}}

\def\ord{{\mathrm{ord}\,}}

\def\({\left(}
\def\){\right)}
\def\fl#1{\left\lfloor#1\right\rfloor}
\def\rf#1{\left\lceil#1\right\rceil}

\def\mand{\qquad \mbox{and} \qquad}

\newcommand{\comm}[1]{\par{\fbox{#1}}}




\hyphenation{re-pub-lished}

\mathsurround=1pt

\def\bfdefault{b}
\overfullrule=5pt

\def \F{{\bbbf}}
\def \K{{\bbbk}}
\def \Z{{\bbbz}}
\def \Q{{\bbbq}}
\def \R{{\bbbr}}
\def \C{{\bbbc}}
\def\Fp{\F_p}
\def \fp{\Fp^*}

\title{On the Number of Solutions of Exponential Congruences}

\author{
{\sc Antal Balog}\\
{Alfr{\'e}d R{\'e}nyi Institute of Mathematics}\\
{Hungarian Academy of Sciences}\\
{H-1364 Budapest, P.O. Box: 127, Hungary}\\
{\tt balog@renyi.hu} \\
\\
{\sc Kevin A. Broughan} \\
{Department of Mathematics}\\
{University of Waikato}\\
{Private Bag 3105, Hamilton, New Zealand}\\
{\tt kab@waikato.ac.nz} \\
\\
{\sc Igor E.~Shparlinski} \\
{Department of Computing}\\
{Macquarie University} \\
{Sydney, NSW 2109, Australia} \\
{\tt igor@comp.mq.edu.au}}

\date{\today }
\maketitle

\begin{abstract}
For a prime $p$ and an integer $a \in \Z$ we obtain
nontrivial upper bounds on
the number of solutions to the congruence
$x^x \equiv a \pmod p$,  $1 \le x \le p-1$.
We use these estimates to estimate the number of solutions
to the congruence
$x^x \equiv y^y \pmod p$,  $1 \le x,y \le p-1$,
which is of cryptographic relevance.
\end{abstract}


 \section{Introduction}

For a prime $p$ and an integer $a \in \Z$ we denote by
$N(p;a)$
the number of solutions to the congruence
\begin{equation}
\label{eq:Cong}
x^x \equiv a \pmod p,  \qquad  1 \le x \le p-1.
\end{equation}
Obviously only the case of $\gcd(a,p)=1$ is of interest.

We note that other than the result Crocker~\cite{Croc} showing that there
are at least $\lfloor \sqrt{(p-1)/2} \rfloor$ incongruent values of
$x^x \pmod p$ when $1\le x\le p-1$ and our estimates,
little appears to be known about the solutions to~\eqref{eq:Cong}.
The function $x \mapsto x^x \pmod p$,
is also used in some cryptographic protocols
(see~\cite[Sections~11.70 and 11.71]{MOV}), so certainly
deserves further investigation, see also~\cite{HoldMor2}
for various conjectures concerning this function.

Here we suggest several approaches to studying this congruence
and derive some upper bounds for $N(p;a)$.

Our first bound is nontrivial if  $a$ is of small multiplicative
order, which in the particular case when $a=1$, 
takes the form $N(p;a)\le p^{1/3+o(1)}$ as $p\to \infty$. The second bound  is
nontrivial if $a$ is of large multiplicative order, which in the
particular case when $a$ is a primitive root modulo $p$, takes
the form
$N(p;a) \le p^{11/12+o(1)}$ as $p\to \infty$.

Furthermore, both bounds combined imply that as $p\to \infty$,
we have the uniform estimate
\begin{equation}
\label{eq:Any a} N(p;a) \le p^{12/13+o(1)}.
\end{equation}

Finally, we estimate the number of solutions $M(p)$ to the
symmetric congruence
\begin{equation}
\label{eq:Cong Sym}
x^x \equiv y^y \pmod p,  \qquad  1 \le x,y \le p-1,
\end{equation}
which has been considered by Holden \& Moree~\cite{HoldMor2}
in their study of short cycles in the iterations of the
discrete logarithm modulo $p$, see also~\cite{Hold,HoldMor1}.
However, no nontrivial estimate of $M(p)$
 has been known prior to
this work. Clearly
\begin{equation}
\label{eq:N and T}
M(p) = \sum_{a=1}^{p-1} N(p;a)^2.
\end{equation}
Thus using the bound~\eqref{eq:Any a} and
the identity
\begin{equation}
\label{eq:Sum T}
 \sum_{a=1}^{p-1} N(p;a) = p-1,
\end{equation}
we immediately derive
\begin{equation}
\label{eq:N easy} M(p) \le p^{25/13+o(1)}.
\end{equation}
However here we obtain a slightly stronger bound, namely
$$
 M(p) \le p^{48/25+o(1)}.
$$

Surprisingly enough, besides elementary number theory arguments,
the bounds derived here rely on some results and arguments from additive
combinatorics, in particular on results of
Garaev~\cite{Gar}.

For an integer $m\ge 1$ we use $\Z_m$ to denote the residue ring
modulo $m$ and we use $\Z_m^*$  to denote the unit group of
$\Z_m$.

Note that without the condition $1\leq x \leq p-1$ (needed in the
cryptographic application) there are always many solutions. Let
$g$ be a primitive root modulo $p$. For any element $a\in\Z_p^*$
(and so for any integer $a \not \equiv  0 \pmod p$) we use $\ind
a$ for its discrete logarithm modulo $p$, that is, the unique
residue class $v\pmod{p-1}$ with
$$
g^v \equiv a \pmod p.
$$
Now, if
for a primitive root $g$ we have
$$
x\equiv p~\ind a-(p-1)g \pmod{p(p-1)},
$$
then
$$
x^x\equiv g^{p~\ind a-(p-1)g}\equiv {(g^p)}^{\ind a}\cdot
(g^{-g})^{p-1}\equiv a\pmod p.
$$

\section{Elements of Small Order}

We need to recall some notions and results from additive combinatorics.

 For a prime $p$ and a set
$\cA \subseteq \Z_p^*$ we define the sets
\begin{equation*}
\cA + \cA  =\{a_1 + a_2\ : \ a_1,   a_2 \in \cA\},   \quad
\cA\cdot\cA  = \{a_1a_2\ : \ a_1, a_2 \in \cA\}.
\end{equation*}

Our bound on $N(p,a)$ makes use of  the
following  estimate of    Garaev~\cite[Theorem~1]{Gar}.

\begin{lemma}
\label{lem:Gar}
For any set $\cA \subseteq \Z_p^*$,
$$
\#(\cA + \cA)  \cdot  \#(\cA \cdot \cA)
\gg \min\left\{p\#\cA, \frac{\(\#\cA\)^4}{p}\right\}.
$$
\end{lemma}

Let $\ord a$ denote the multiplicative order of $a \in \Z_p^*$.

\begin{theorem}
\label{thm:Small} Uniformly over $t \mid p-1$, we have,
as $p\to \infty$,
$$\sum_{\substack{a \in \Z_p^*\\ \ord a\mid t}}N(p;a)  \le
\max\{t ,  p^{1/2}t^{1/4}\} p^{o(1)}.
$$
\end{theorem}

\begin{proof}  Fix a primitive root $g$ mod $p$.
 The union of non-zero residue classes $a$ with $\ord a\mid t$ of all
the solutions to~\eqref{eq:Cong} is precisely the set of solutions to
\begin{equation}
\label{eq:tCong} x^{tx} \equiv 1 \pmod p,  \qquad  1 \le x \le
p-1.
\end{equation}
This congruence is equivalent to
$$tx ~\ind x \equiv 0 \pmod {p-1},$$
or if we put
$$
T = \frac{p-1}t
$$
to
$$x ~\ind x \equiv 0 \pmod{T},$$
or after fixing $d \mid T$ and considering only the
solutions to~\eqref{eq:tCong} with
$$
\gcd(x,T) = d,
$$
they can be written as $x=dy$ and satisfy
\begin{equation}
\label{eq:tyCong} \ind (dy)\equiv 0 \pmod{T_d}, \qquad
1\leq y\leq D,\qquad \gcd\(y,T_d\)=1.
\end{equation}
where
$$
T_d= \frac{T}{d} \mand  D = \frac{p-1}d.
$$
Let us denote by $\cY_d$ the set of integers $y$
satisfying~\eqref{eq:tyCong}, and by $\cW_d$ the set of the residue classes
mod $p$ represented by the elements of $\cY_d$. Obviously
$\#\cY_d=\#\cW_d$, and we have
\begin{equation}
\label{eq:Basic}\sum_{\substack{a \in \Z_p^*\\
\ord a\mid t}}N(p;a)  = \sum_{d\mid T}\#\cY_d =
\sum_{d\mid T}\#\cW_d.
\end{equation}

First note that
\begin{equation}
\label{eq:SumSet} \# \(\cW_{d} + \cW_{d}\) \le \# \(\cY_{d} +
\cY_{d}\) \le 2D
\end{equation}
from the second condition in~\eqref{eq:tyCong}.

Furthermore, the product set of $\cW_d$ is contained in
$$\{w\in\Z_p^*:\,\ind(d^2w)\equiv 0 \pmod{T_d}\},$$ and so
\begin{equation}
\label{eq:ProdSet} \# \(\cW_{d} \cdot \cW_{d}\)  \le \frac{p-1}{T_d} = dt.
\end{equation}

Hence, applying Lemma~\ref{lem:Gar}
and using the bounds~\eqref{eq:SumSet} and~\eqref{eq:ProdSet}
we see that
$$
\min\left\{p\#\cW_{d}, \frac{\(\#\cW_{d}\)^4}{p}\right\} \ll pt.
$$
Hence
\begin{equation}
\label{eq:Bound}
\# \cW_{d}\ll  \max\{t ,  p^{1/2}t^{1/4}\} .
\end{equation}

Recalling the bound on the divisor function  $\tau(k)$
\begin{equation}
\label{eq:Div Bound}
\tau(k) = \sum_{d \mid k} 1 = k^{o(1)},
\end{equation}
see~\cite[Theorem~315]{HardyWright}, and using~\eqref{eq:Bound}
in~\eqref{eq:Basic}, we conclude the proof.
\end{proof}

\begin{cor}
\label{cor:Small} Uniformly over $t \mid p-1$ and
all integers $a$
with $\gcd(a,p)=1$  of multiplicative order $\ord a  = t$,
we have, as $p\to \infty$,
$$N(p;a)  \le  \max\{t ,  p^{1/2}t^{1/4}\} p^{o(1)}.
$$
\end{cor}


Next we show that if $t$ is very small then the bound of
Theorem~\ref{thm:Small} can be improved. For example, this applies
to the most interesting special case of the
congruence~\eqref{eq:Cong}, namely the case $a=1$.

\begin{theorem}
\label{thm:Very Small} Uniformly over $t \mid p-1$, we have,
as $p\to \infty$,
$$\sum_{\substack{a \in \Z_p^*\\ \ord a\mid t}}N(p;a)
\le  p^{1/3 + o(1)} t^{2/3}.
$$
\end{theorem}

\begin{proof} We follow
the proof of Theorem~\ref{thm:Small} up to~\eqref{eq:ProdSet}, but
finish the argument in a different way to derive a new bound for
$\#\cY_d$. Let us define
$$
s(b)=\#\{(y_1,y_2)\ : \ y_1,y_2\in \cY_d, \  y_1y_2 \equiv b \pmod
p\}.
$$
First note that $s(b)>0$ only when $b\in \cW_d\cdot\cW_d$, and so
\begin{equation}
\label{eq:Prelim} (\# \cY_d)^2=\sum_{b\in \Z_p} s(b) \le
\#\(\cW_d \cdot \cW_d\) \max_{b\in\Z_p} s(b).
 \end{equation}

If $(y_1,y_2)$ is counted in $s(b)$ then on the one hand $y_1y_2
\equiv b \pmod p$, on the other hand $1\leq
y_1y_2\leq D^2$ (where as before $D = (p-1)/d$),
therefore $y_1y_2=b+kp$, where $0\leq k< \frac{p}{d^2}$. Thus the product $y_1y_2$ can take at most
$p/d^2 + 1$ possible values $y_1y_2=z$ and once $z$ is fixed,
there are $\tau(z) = z^{o(1)} = p^{o(1)}$ possibilities for the
pair $(y_1,y_2)$, see~\eqref{eq:Div Bound}. Thus
$$
s(b) \le  (p/d^2 + 1)  p^{o(1)},
$$
which after inserting in~\eqref{eq:Prelim}
and recalling~\eqref{eq:ProdSet} yields
\begin{equation}
\label{eq:Bound 3} \# \cY_d \le \((pt/d)^{1/2} + (td)^{1/2}\)
p^{o(1)}.
\end{equation}

For $d \le p^{1/3} t^{-1/3}$ we use $\#\cY_d\leq dt$ from the
first condition of~\eqref{eq:tyCong} and for $d \ge p^{2/3}
t^{-1/3}$ we use $\#\cY_d\leq D$ from the second
condition of~\eqref{eq:tyCong}. Therefore we obtain
$$
\# \cY_d \ll p^{1/3} t^{2/3} \mand \#\cY_d \ll p^{1/3} t^{1/3},
$$
respectively.

Finally, for $  p^{1/3} t^{-1/3} \le d \le p^{2/3} t^{-1/3}$ we
use~\eqref{eq:Bound 3} to derive
$$
\# \cY_d  \le \(p^{1/3} t^{2/3} + p^{1/3} t^{1/3}\)p^{o(1)} =
p^{1/3+o(1)} t^{2/3}.
$$

Using these bounds with~\eqref{eq:Div Bound} in~\eqref{eq:Basic}
we conclude the proof.
\end{proof}

\begin{cor}
\label{cor:VerySmall}  Uniformly over $t\mid p-1$ and  all
integers $a$ with $\gcd(a,p)=1$  of multiplicative order $\ord a
= t$, we have, as $p\to \infty$,
$$N(p;a)  \le  p^{1/3 + o(1)} t^{2/3}.
$$
\end{cor}

\section{Elements of Large Order}

Here we use a different argument, which is similar to the one
used in~\cite{BouShp}, and a bound of~\cite{CFKLLS}, on the number
of solutions of an exponential congruence, plays the crucial role.
However, this approach is effective only for values of $a$ of
sufficiently large order.

We recall the following estimate, given in~\cite[Lemma~7]{CFKLLS},
on the number of zeros
of sparse polynomials over a finite field $\F_q$ of $q$
elements.

\begin{lemma}
\label{lem:SprPol}
For
$n \ge 2$ given
elements  $a_1, \ldots, a_n \in \F_q^*$ and
integers $k_1, \ldots, k_n$ in $\Z$
let us denote by $Q$ the number of solutions of the equation
$$
\sum_{i=1}^n a_i X^{k_i} = 0, \qquad X \in \F_q^*.
$$
Then
$$
Q \le 2 q^{1 - 1/(n-1)} \Delta^{1/(n-1)}
+ O\(q^{1 - 2/(n-1)} \Delta^{2/(n-1)} \),
$$
where
$$
\Delta =  \min_{1 \le i \le n} \max_{j \ne i} \gcd(k_j - k_i, q-1).
$$
\end{lemma}

We are now ready  to
prove the main result of this section.

\begin{theorem}
\label{thm:Large}
 Uniformly over $t \mid p-1$ and all integers $a$
with $\gcd(a,p)=1$  of multiplicative order $\ord a  = t$,
we have, as $p\to \infty$,
$$
N(p;a)  \le  p^{1 + o(1)}t^{-1/12}.
$$
\end{theorem}

\begin{proof} Let $a$ be a non-zero residue class modulo $p$ of
multiplicative order $t\mid p-1$.
As before, we put
$$
T = \frac{p-1}t
$$

Clearly, there is a primitive root $g$
modulo $p$ with $a\equiv g^{T}\pmod p$. Using the discrete
logarithm to base $g$, the congruence~\eqref{eq:Cong} is equivalent to
$$
x ~\ind x \equiv  T \pmod {p-1}.
$$
Note the condition $\gcd(x,p-1) \mid T$. After fixing
$d \mid T$ and considering only the solutions
to~\eqref{eq:Cong} with $ \gcd(x,p-1) = d$, they can be written
as $x=dy$ and satisfy
$$
 y~\ind (dy)\equiv T_d
\pmod{D}, \quad 1\leq y\le D, \quad
\gcd(y,D)=1,
$$
where, as before,
$$
T_d= \frac{T}{d} \mand  D = \frac{p-1}d.
$$

Note that $t\mid D$. The congruence
$yz\equiv 1 \pmod{D}$ defines a
one--to--one correspondence between the integers $\{1\leq y\leq
D\ : \ \gcd(y,D)=1\}$ and $z\in\Z_D^*$.

Furthermore, the relation $yz\equiv 1 \pmod D$ defines
a one--to--$M_d$
correspondence between the set
 $\{1\leq y\leq
D\ : \ \gcd(y,D)=1\}$ and $z\in\Z_{p-1}^*$, where $M_d$ is the number of residue classes in
$\Z_{p-1}^*$ in the form $z+kD$. These residue classes are
automatically coprime to $D$, but we have to ensure  that
they are coprime to $d$ as well (and thus belong to $\Z_{p-1}^*$).
Thus using $\mu(k)$ to denote
the M{\"o}bius function, by~\cite[Theorem~263]{HardyWright}
(which is essentially the inclusion-exclusion principle)
we obtain
\begin{eqnarray*}
M_d&=&\sum_{k=1}^d\,\sum_{f\mid\gcd(z+kD,d)}\mu(f) =
\sum_{f\mid d}\mu(f)\sum_{\substack{k=1\\z+kD\equiv 0\pmod f}}^d 1\\
&=& \sum_{\substack{f\mid d\\\gcd(f,D)=1}}\mu(f)\frac df =
d\,\frac{\varphi(m)}m,
\end{eqnarray*}
where $\varphi(k)$ is the Euler function and
$m$ is the product of primes $q$ with $q\mid d$ and $q\nmid D$,
see~\cite[Equation~(16.3.1)]{HardyWright}.
In particular $m\le d\le p$ and recalling the well-known estimate
on the Euler function, see~\cite[Theorem~328]{HardyWright} we obtain
$$
M_d = d p^{o(1)}.
$$

{From} now on the integer $1\le y\le D$ and the residue class
$z\in \Z_{p-1}^*$ with or without subscripts are always
connected by $yz\equiv 1\pmod D$, even if this is not explicitly
stated.

Let us define
$$
\cZ_d=\{z\in\Z_{p-1}^*\ :\ \ind (dy)\equiv Dz/t\pmod D,\ 1\leq y\leq D\}.
$$
(we recall our convention that we always have $yz\equiv 1 \pmod D$).
We have
\begin{equation}\label{eq:zBasic}
N(p,a)=\sum_{d\mid T}\frac {1}{M_d}\#\cZ_d\leq
p^{o(1)}\sum_{d\mid T}\frac 1d\#\cZ_d.
\end{equation}

The congruence $\ind (dy)\equiv Dz/t \pmod D$ is equivalent to
$$
dy\equiv \rho g^{Dz/t}\pmod p,
$$
for some $\rho\in\Z_p^*$ with $\rho^d\equiv 1\pmod p$. Thus we split
$\cZ_d$ into subsets $\cZ_{d,\rho}$ getting
\begin{equation}\label{eq:zroBasic}
\#\cZ_d=\sum_{\rho^d\equiv 1\pmod p}\#\cZ_{d,\rho},
\end{equation}
where
$$\cZ_{d,\rho}=\{z\in\Z_{p-1}^*\ :\ dy\equiv\rho g^{Dz/t} \pmod p,\  1\leq y\leq D\}
$$
(and again we recall our convention that  $yz\equiv 1 \pmod D$).


Clearly,
$$
(\#\cZ_{d,\rho})^2=\#\{z_1,z_2\in\Z_{p-1}^*\ :\ dy_j\equiv\rho
g^{Dz_j/t}\pmod p,\ j=1,2\}.
$$
We have by adding the two congruences  that
\begin{eqnarray*}
\lefteqn{(\#\cZ_{d,\rho})^2}\\
& & \quad \le \#\{z_1,z_2\in\Z_{p-1}^*\ :\ d(y_1+y_2)\equiv \rho
\(g^{Dz_1/t}+g^{Dz_2/t}\)\pmod p\}\\
& &  \quad =  \sum_{v\in \Z}\#\{z_1,z_2\in\Z_{p-1}^*\ :\ d(y_1+y_2)=v,\\
& & \qquad \qquad\qquad \qquad\qquad \qquad  \quad \rho
\(g^{Dz_1/t}+g^{Dz_2/t}\)\equiv v\pmod p\}.
\end{eqnarray*}

The sum over $v \in \Z$ is empty
unless $v=dw$, where $2\leq w\leq 2D$ and we get
by the Cauchy--Schwarz inequality that
\begin{eqnarray*}
\lefteqn{(\#\cZ_{d,\rho})^4 \le 2D\#\{z_1, z_2, z_3, z_4\in\Z_{p-1}^*\
: \ d(y_1+y_2)=d(y_3+y_4)}\\
& & \qquad \qquad  \qquad \equiv \rho
\(g^{Dz_1/t}+g^{Dz_2/t}\)\equiv \rho
\(g^{Dz_3/t}+g^{Dz_4/t}\)\pmod p\}.
\end{eqnarray*}

Clearly, when $z_1, z_2, z_3, z_4\in\Z_{p-1}^*$ are fixed,
then the condition
\begin{eqnarray*}
\lefteqn{
 d(y_1+y_2)=d(y_3+y_4) }\\
& & \qquad \qquad  \qquad
\equiv \rho
\(g^{Dz_1/t}+g^{Dz_2/t}\)\equiv \rho
\(g^{Dz_3/t}+g^{Dz_4/t}\) \pmod p
\end{eqnarray*}
defines $\rho$ uniquely.  Hence
\begin{eqnarray*}
\lefteqn{
\sum_{\rho^d\equiv 1\pmod p}(\#\cZ_{d,\rho})^4}\\
& & \qquad  \le 2D\#\{z_1, z_2, z_3, z_4\in\Z_{p-1}^*\
: \ y_1+y_2=y_3+y_4, \\
& & \qquad \qquad\qquad \qquad\qquad
g^{Dz_1/t}+g^{Dz_2/t} \equiv  g^{Dz_3/t}+g^{Dz_4/t}\pmod p\}.
\end{eqnarray*}
Relaxing the condition $y_1+y_2=y_3+y_4$ to $y_1+y_2\equiv
y_3+y_4\pmod D$ only increases the number
of solution (but allows us to think about $y_j$ as a residue class
modulo $D$ defined by $y_jz_j\equiv 1 \pmod D$, $j=1,2,3,4$).
Thus
\begin{eqnarray*}
\lefteqn{
\sum_{\rho^d\equiv 1\pmod p}(\#\cZ_{d,\rho})^4}\\
& & \qquad  \le 2D\#\{z_1, z_2, z_3, z_4\in\Z_{p-1}^*\
: \ y_1+y_2\equiv y_3+y_4\pmod D, \\
& & \qquad \qquad\qquad \qquad\qquad  \quad
g^{Dz_1/t}+g^{Dz_2/t} \equiv  g^{Dz_3/t}+g^{Dz_4/t}\pmod p\}.
\end{eqnarray*}

Finally, after the substitution $z_j \to wz_j$ for $w \in
\Z_{p-1}^*$ (and thus $y_j \to w^{-1}y_j$),  $j=1,2,3,4$,
where $w^{-1}$ is defined
modulo $D$, we obtain that any solution is computed with $\varphi(p-1)$
multiplicity, that is
\begin{equation}\label{eq:U W}
\begin{split}
&\sum_{\rho^d\equiv 1\pmod p}\, (\#\cZ_{d,\rho})^4
\le
\frac{2D}{\varphi(p-1)} \,
\#\{z_1, z_2, z_3, z_4, w \in\Z_{p-1}^*\ :\\
&  \qquad \qquad y_1+y_2\equiv y_3+y_4 \pmod D,\\
 &  \qquad \qquad \quad  (g^w)^{Dz_1/t}+(g^w)^{Dz_2/t} \equiv  (g^w)^{Dz_3/t}+(g^w)^{Dz_4/t}
\pmod p\}.
\end{split}
\end{equation}
Writing $X\equiv g^w \pmod p$ and
$k_j=Dz_j/t=(p-1)z_j/dt = T_d z_j$, after
fixing $z_1, z_2, z_3, z_4$,
the number of $w\in \Z_{p-1}^*$ satisfying the
congruence in~\eqref{eq:U W} is bounded by the number of solutions
to the congruence $X^{k_1}+X^{k_2}\equiv X^{k_3}+X^{k_4}\pmod p$,
and this is bounded in Lemma~\ref{lem:SprPol}, applied with $n=4$,
by $O\(p^{2/3} \Delta^{1/3}\)$, where
\begin{eqnarray*}
\Delta &=& \min_{1 \le i < j \le 4}
\gcd\(T_d(z_i - z_j), p-1\)= T_d \min_{1 \le i < j\le 4}
\gcd\(z_i - z_j, dt\).
\end{eqnarray*}

For every fixed $i,j$, $1 \le i < j \le 4$
and $\delta \mid dt$ there are $(p-1)^2/\delta$
choices for $(z_i,z_j)$ with
$$\gcd(z_i - z_j, dt) = \delta.
$$
When
$z_i$ and $z_j$ are fixed the congruence $y_1 + y_2 \equiv y_3 +
y_4\pmod D$ implies that there are $dp^{1+o(1)}$ choices for
the remaining two variables. (Recall that each $y$ determines
$M_d=dp^{o(1)}$ different choices of $z$.) Thus, putting
everything together in~\eqref{eq:U W}
and recalling~\eqref{eq:Div Bound}, we obtain
\begin{eqnarray*}
\lefteqn{\sum_{\rho^d\equiv 1\pmod p}\,(\#\cZ_{d,\rho})^4 \le \frac{2D}{\varphi(p-1)} \sum_{\delta\mid
dt}p^{2/3} (T_d \delta)^{1/3}  \frac{(p-1)^2}{\delta} dp^{1+o(1)}}  \\
& & \qquad \qquad \quad  = d Dp^{8/3+o(1)} T_d^{1/3}\sum_{\delta\mid
dt}\delta^{-2/3}  = p^{11/3+o(1)} T_d^{1/3}
= \frac{p^{4+o(1)}}{(dt)^{1/3}}.
\end{eqnarray*}
Putting this to~\eqref{eq:zroBasic}, we get by
the H\"older inequality
$$
\#\cZ_d\le d^{3/4}\left(\sum_{\rho^d\equiv 1\pmod
p}\,(\#\cZ_{d,\rho})^4\right)^{1/4}\le
\frac{p^{1+o(1)}}{t^{1/12}}\,d^{2/3}.
$$
Finally~\eqref{eq:zBasic} and~\eqref{eq:Div Bound} gives
$$
N(p,a)\le \sum_{d\mid
(p-1)/t}\,\frac{p^{1+o(1)}}{t^{1/12}d^{1/3}}\le
\frac{p^{1+o(1)}}{t^{1/12}},
$$
and we conclude the proof.
\end{proof}

\section{Symmetric Congruence}

We now improve the bound~\eqref{eq:N easy}
on the number of solutions to the symmetric congruence~\eqref{eq:Cong Sym}.

\begin{theorem}
\label{thm:Sym Cong}
We  have, as $p\to \infty$.
$$
M(p)  \le  p^{48/25 + o(1)}.
$$
\end{theorem}

\begin{proof}
{From}~\eqref{eq:N and T} we obtain
$$
M(p) \le \sum_{t \mid p-1}
\sum_{\substack{a \in \Z_p^*\\ \ord a=t}}N(p;a)^2.
$$
We fix some parameter $\vartheta$ and for $t \le \vartheta$
we use Theorem~\ref{thm:Small} to estimate
\begin{eqnarray*}
\sum_{\substack{a \in \Z_p^*\\ \ord a=t}}N(p;a)^2
& \le & \(\sum_{\substack{a \in \Z_p^*\\ \ord a=t}}N(p;a)\)^2\\
& \le &\max\{t^2 p^{o(1)}, p^{1+ o(1)} t^{1/2}\} \le\max\{\vartheta^2
p^{o(1)}, p^{1+ o(1)} \vartheta^{1/2}\}.
\end{eqnarray*}

For  $t \ge \vartheta$  we use Theorem~\ref{thm:Large}
together with~\eqref{eq:Sum T} to estimate
$$
\sum_{\substack{a \in \Z_p^*\\ \ord a=t}}N(p;a)^2
\le  p^{1 + o(1)}t^{-1/12}
\sum_{\substack{a \in \Z_p^*\\ \ord a=t}}N(p;a)
\le  p^{2 + o(1)} \vartheta^{-1/12}.
$$

Taking
$$
\vartheta = p^{24/25}
$$
to balance the above estimates,
we obtain the bound
$$
\sum_{\substack{a \in \Z_p^*\\ \ord a=t}}N(p;a)^2 \le  p^{48/25 +
o(1)}
$$
and using~\eqref{eq:Div Bound}, we  conclude the proof.
\end{proof}

 \section{Concluding Remarks}

Clearly Theorem~\ref{thm:Small} is nontrivial provided that $t\le
p^{1-\varepsilon}$ for some $\varepsilon >0$, while
Theorem~\ref{thm:Large} is nontrivial provided $t\ge
p^{\varepsilon}$, for an arbitrary $\varepsilon>0$ and a
sufficiently large $p$.  In particular, using
Corollary~\ref{cor:Small} for  $t\le p^{12/13}$ and
Theorem~\ref{thm:Large} for $t> p^{12/13}$, we
derive~\eqref{eq:Any a}.

It is also easy to see that all but $o(p)$ elements $a \in \Z_p^*$ are
of multiplicative order $t= p^{1+o(1)}$. Thus for almost all  $a \in \Z_p^*$
we have $N(p;a)  \le p^{11/12+o(1)}$ by Theorem~\ref{thm:Large}.

Similar results can also be established for several other congruences.
For example, the same arguments as those used in the proof of
Theorem~\ref{thm:Very Small} imply that the congruence
$$
x^{x-1} \equiv 1 \pmod p,  \qquad  1 \le x \le p-1,
$$
has $O\(p^{1/3 + o(1)}\)$ solutions. This means that the function
$x \mapsto x^x \pmod p$ has $O(p^{1/3 + o(1)})$ fixed points in
the interval $1 \le x \le p-1$.

\section*{Acknowledgements}

Research of A.~B. was supported in part by Hungarian National
Science Foundation Grants K72731 and K81658
and that of I.~S.\ was supported in part by Australia Research
Council Grants DP0556431 and DP0881473.

\end{document}